\documentclass{svjour3}                     
\smartqed  
\usepackage{graphicx}
\usepackage{hyperref}
\usepackage{mathptmx}      
\usepackage{amsmath,amssymb,latexsym}
\usepackage{color}
\newcommand{\NN}{\mathbb{N}}
\newcommand{\RR}{\mathbb{R}}

\newcommand{\cS}{\mathcal S}

\newcommand{\cP}{\mathcal P}

\newcommand{\raw}{\rightarrow}
\journalname{Foundations of Computational Mathematics}
\begin{document}

\title{Nodal bases for the serendipity family of finite elements
\thanks{The published version of this pre-print appears in The Journal of the Society for the Foundations of Computational Mathematics.
The final publication is available at Springer via:\\
\text{}\qquad\qquad\qquad\qquad\qquad\href{http://dx.doi.org/10.1007/s10208-016-9305-0}{\texttt{http://dx.doi.org/10.1007/s10208-016-9305-0}}
}
\thanks{AG was supported in part by NSF Award 1522289.}
}



\author{Michael S. Floater  \and
       Andrew Gillette
}


\institute{M. Floater \at
  Department of Mathematics, 
  University of Oslo, 
  PO Box 1053, Blindern, 
  0316 Oslo, Norway, 
  michaelf@math.uio.no
           \and
  A. Gillette\at
  Department of Mathematics, 
  University of Arizona, 
  617 N. Santa Rita Avenue,
  Tucson, AZ, USA 85721,
  agillette@math.arizona.edu
}

\date{}

\maketitle

\begin{abstract}
Using the notion of multivariate lower set interpolation,
we construct nodal basis functions for the serendipity
family of finite elements, of any order and any dimension.
For the purpose of computation, we also show how to express these functions
as linear combinations of tensor-product polynomials.
%
%
\end{abstract}


\section{Introduction}\label{sec:intro}

The serendipity family of~$C^0$ finite elements is commonly used on
cubical and parallelepiped meshes in two and three dimensions as a means to reduce the computational effort required by tensor-product elements.
The number of basis functions of a tensor-product element of order~$r$ in~$n$ dimensions is~$(r+1)^n$, while for a serendipity element it is asymptotically~$\sim r^n / n!$ for large~$r$, which represents a reduction of 50\% in 2-D and 83\% in 3-D.
In this paper, we construct basis 
functions for serendipity elements of any order $r\geq 1$ in any number 
of dimensions $n\geq 1$, that are interpolatory at specified nodes and 
can be written as linear combinations of tensor-product 
polynomials (see equation (\ref{eq:nodalformula})).
The benefits and novelty of our approach are summarized as follows:
\begin{itemize}
\item \textbf{Flexible node positioning.}
Our approach constructs nodal basis functions for any arrangement of points on the $n$-cube that respects the requisite association of degrees of freedom with sub-faces.
In particular, we allow a symmetric arrangement of points that remains invariant under the symmetries of the $n$-cube. 
\item \textbf{Tensor product decomposition.} 
The basis functions we define are can be written as linear combinations of standard tensor product basis functions, with coefficients prescribed by a simple formula based on the geometry of a lower set of points associated to superlinear monomials.
\item \textbf{Dimensional nesting.} 
The restriction of our basis functions for a $n$-cube to one of its $s$-dimensional faces coincides with the definition of our basis functions for an $s$-cube.
\end{itemize}

Serendipity elements have appeared in various
mathematical and engineering texts, typically for
small $n$ such as $n=2$ and $n=3$,
and for small $r$; see~\cite{BS2002,Ci02,G,H1987,M1990,SB1991,SF73,LZ,Z1998}.
A common choice for the basis functions is a nodal (Lagrange) basis,
which is an approach that has also been studied in the
approximation theory literature.
For example, Delvos \cite{D} applied his `Boolean
interpolation' to construct a nodal basis for the case $n=3$ and $r=4$.
Other bases have been considered, such as products
of univariate Legendre polynomials, as in the work of
Szab{\'o} and Babu{\v{s}}ka~\cite{SB1991}.

It was relatively recently that the serendipity spaces
were chacterized precisely
for arbitrary $n$ and $r$, by Arnold and Awanou~\cite{AA}.
They derived the polynomial space and its dimension, and
also constructed a unisolvent set of degrees of freedom
to determine an element uniquely.
For the $n$-dimensional cube $I^n$, with $I = [-1,1]$, they
defined the serendipity space $\cS_r(I^n)$ as the linear space
of $n$-variate polynomials whose \textit{superlinear degree} is at most $r$.
The superlinear degree of a monomial is its total degree, less the number of variables appearing only linearly in the monomial.
For a face $f$ of $I^n$ of dimension $d \ge 1$,
the degrees of freedom proposed in \cite{AA}
for a scalar function $u$ are of the form
\begin{equation}
\label{eq:srdp-dofs}
u\longmapsto\int_f uq,
\end{equation}
for $q$ among some basis of $\cP_{r-2d}(f)$.
Here, $\cP_s(f)$ is the space of restrictions to $f$
of $\cP_s(I^n)$, the space of $n$-variate polynomials of
degree $\le s$.
These degrees of freedom were shown to be unisolvent
by a hierarchical approach through the $n$ dimensions:
the degrees of freedom at the vertices of $I^n$ are determined
first (by evaluation);
then the degrees of freedom on the 1-dimensional faces (edges),
then those on the 2-dimensional faces, etc., finishing with
those in the interior of~$I^n$.

The approach of \cite{AA} has the advantage that the degrees of freedom
on any face $f$ of any dimension $d$
can be chosen independently of those on another face, of the
same or of different dimension.
Implementing a finite element method using these degrees of freedom, however, requires a set of `local basis functions' that are associated to the integral degrees of freedom in some standardized fashion.
The lack of simple nodal basis functions for this purpose has limited the broader use and awareness of serendipity elements.

The purpose of this paper is to show that by applying
the notion of lower set interpolation in approximation theory
and choosing an appropriate Cartesian grid in $I^n$,
a nodal basis can indeed be constructed for the serendipity space
$\cS_r(I^n)$ for any $n$ and $r$. The interpolation nodes
are a subset of the points in the grid.
The restrictions of the basis functions to any
$d$-dimensional face are themselves basis functions of the same type
for a $d$-cube, yielding $C^0$ continuity between adjacent elements.

If we keep all the nodes distinct, it is not possible
to arrange them in a completely symmetric way for $r \ge 5$.
However, lower set interpolation also applies to derivative data,
and by collapsing interior grid coordinates to the midpoint of $I$,
we obtain a Hermite-type basis of functions
that are determined purely by symmetric interpolation conditions
for all $n$ and $r$.

Lower set interpolation can also be expressed as
a linear combination of tensor-product interpolants
on rectangular subsets of the nodes involved \cite{DF}.
We derive an explicit formula for the coefficients in the
serendipity case, which could be used for evaluation of
the basis functions and their derivatives.

\section{Interpolation on lower sets}\label{sec:lower}

A multi-index of $n$ non-negative integers will be denoted by
\[ \alpha = (\alpha_1,\alpha_2,\ldots,\alpha_n) \in \NN_0^n. \]
For each $j=1,\ldots,n$, choose grid coordinates
$x_{j,k}\in\RR$ for all $k \in \NN_0$,
not necessarily distinct. These coordinates determine
the grid points
\begin{equation}
\label{eq:x-alpha}
x_\alpha := (x_{1,\alpha_1}, x_{2,\alpha_2}, \ldots, x_{n,\alpha_n})
  \in \RR^n, \qquad \alpha \in \NN_0^n.
\end{equation}
The \emph{left multiplicity} of $\alpha \in \NN_0^n$ with respect
to the $x_{j,k}$ is defined to be the multi-index
\[ \rho(\alpha) := (\rho_1(\alpha),\ldots,\rho_n(\alpha)) \in \NN_0^n, \]
where
\begin{equation}
\label{eq:rhoj-def}
\rho_j(\alpha) := \# \{k < \alpha_j : x_{j,k} = x_{j,\alpha_j} \}.
\end{equation}
Thus $\rho_j(\alpha)$ is the number of coordinates in the sequence
$x_{j,0}, x_{j,1},\ldots, x_{j,\alpha_j-1}$ that are equal to $x_{j,\alpha_j}$.
For each $\alpha\in\NN_0^n$, we associate a linear functional
$\lambda_\alpha$ as follows.
Given any $u:\RR^n\raw\RR$, defined with sufficiently many derivatives
in a neighborhood of $x_\alpha$, let
\[ \lambda_\alpha u := D^{\rho(\alpha)} u(x_\alpha). \]
We call a finite set $L \subset \NN_0^n$ a \emph{lower set}
if $\alpha \in L$ and
$\mu \leq \alpha$ imply $\mu \in L$. 
The partial ordering $\mu \leq \alpha$ means $\mu_j \leq \alpha_j$
for all $j=1,\ldots,d$.
We associate with $L$ the linear space of polynomials
\begin{equation}
\label{eq:PofL}
P(L)= {\rm span}\{x^\alpha : \alpha \in L \},
\end{equation}
where
\begin{equation}\label{eq:monomial}
 x^\alpha := x_1^{\alpha_1} \cdots x_n^{\alpha_n},
\end{equation}
for any point $x = (x_1,\ldots,x_n) \in \RR^n$.

Polynomial interpolation on lower sets
has been studied in
\cite{CCS,BR,D,DF,GS2,Kuntzmann,LL,Muhlbach,Sauer,Werner}
and the following theorem has been established in various special cases
by several authors.

\begin{theorem}
\label{thm:NewtonH}
For any lower set $L \subset \NN_0^n$ and
a sufficiently smooth function $u:\RR^n\raw\RR$,
there is a unique polynomial $p \in P(L)$ that interpolates $u$ in the
sense that
\begin{equation}\label{eq:agrees}
 \lambda_\alpha p = \lambda_\alpha u, \qquad \alpha \in L.
\end{equation}
\end{theorem}
The theorem leads to a basis of $P(L)$ with the basis function 
$\phi_\alpha\in P(L)$, $\alpha\in L$, defined by
\begin{equation}\label{eq:basis}
 \lambda_{\alpha'} \phi_\alpha = \delta_{\alpha, \alpha'}, \qquad \alpha \in L,
\end{equation}
where $\delta_{\alpha,\alpha'}$ is 1 if $\alpha=\alpha'$ and 0 otherwise.
We can then express $p$ as
\[ p(x) = \sum_{\alpha \in L} \phi_\alpha(x) \lambda_\alpha u. \]

\section{Serendipity spaces}\label{sec:seren}

The serendipity space $\cS_r(I^n)$
can be described and partitioned using the
language of lower sets.
The standard norm for a multi-index $\alpha \in \NN_0^n$ is
\[ |\alpha| := \sum_{j=1}^n \alpha_j, \]
which is the degree of the monomial $x^\alpha$ in
(\ref{eq:monomial}).
We will define the \emph{superlinear} norm of $\alpha$ to be
\[ |\alpha|' := \sum_{\substack{j=1 \\ \alpha_j \ge 2}}^n \alpha_j, \]
which is the `superlinear' degree of $x^\alpha$ from~\cite{AA}.
Using this we define, for any $r \ge 1$,
\begin{equation}
\label{eq:Srdef}
S_r := \{\alpha \in \NN_0^n : |\alpha|' \leq r \}.
\end{equation}
Observe that $S_r$ is a lower set since $|\alpha|' \leq |\beta|'$
whenever $\alpha \leq \beta$.
Recalling (\ref{eq:PofL}), we let $\cS_r = P(S_r)$, which
coincides with the definition of $\cS_r$ in~\cite{AA}.

We now partition $S_r$, and hence $\cS_r(I^n)$,
with respect to the faces of $I^n$.
We can index these faces using a multi-index
$\beta \in \{0,1,2\}^n$. For each such $\beta$ we define the face
\[ f_\beta = 
   I_{1,\beta_1} \times I_{2,\beta_2} \times \cdots \times I_{n,\beta_n}, \]
where
\begin{equation}
\label{eq:ijbetaj}
I_{j,\beta_j} := \begin{cases}
                   -1, & \beta_j = 0; \\
                   1, & \beta_j = 1; \\
                   (-1,1), & \beta_j = 2.
                 \end{cases}
\end{equation}
Since $I$ can be written as the disjoint union 
$ I = \{-1\} \cup \{1\} \cup (-1,1)$, we see that $I^n$ can be written as the disjoint union
$$ I^n = \bigcup_{\beta \in \{0,1,2\}^n} f_\beta. $$
Hence, there are $3^n$ faces of all dimensions.
The dimension of the face $f_\beta$ is 
\[ \dim f_\beta = \# \{j : \beta_j = 2\}, \]
and the number of faces of dimension $d$ is
\begin{equation}\label{eq:numfaces}
\#\{f_\beta\subseteq I^n : \dim f_\beta=d\} =  2^{n-d} \binom{n}{d}.
\end{equation}
The $2^n$ vertices of $I^n$ correspond to $\beta\in\{0,1\}^n$,
the $2^{n-1} n$ edges correspond to $\beta$ with exactly one entry equal
to 2, and so forth, up to the single $n$-face, $f_{(2,2,\ldots,2)}$,
the interior of $I^n$.
To partition $S_r$ according to these faces, write $S_r$ as the
disjoint union
\begin{equation}
\label{eq:Sr-part}
S_r = \bigcup_{\beta \in \{0,1,2\}^n} S_{r,\beta},
\end{equation}
where
\[ S_{r,\beta} = \{\alpha \in S_r: \min(\alpha_j,2) 
 = \beta_j,~\text{for}~j=1,\ldots,n\}. \]

\begin{figure}
\centering{
\begin{tabular}{ccc}
\parbox{.17\textwidth}{\includegraphics[width=.20\textwidth]{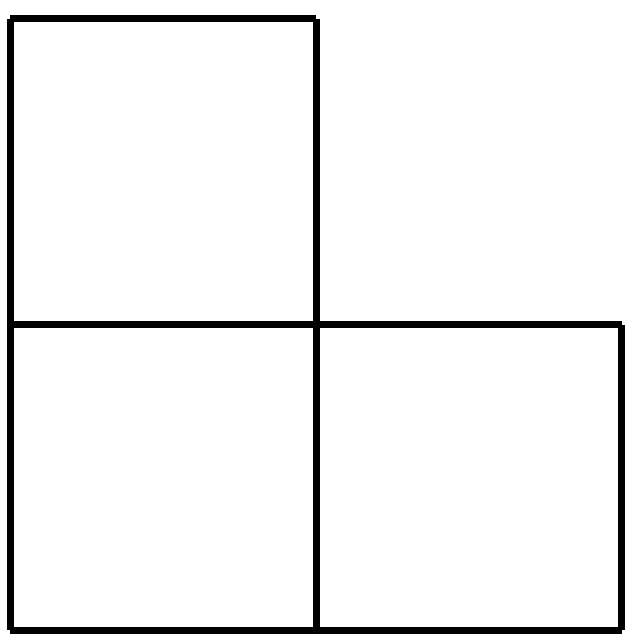}} &
\quad
\parbox{.17\textwidth}{\includegraphics[width=.20\textwidth]{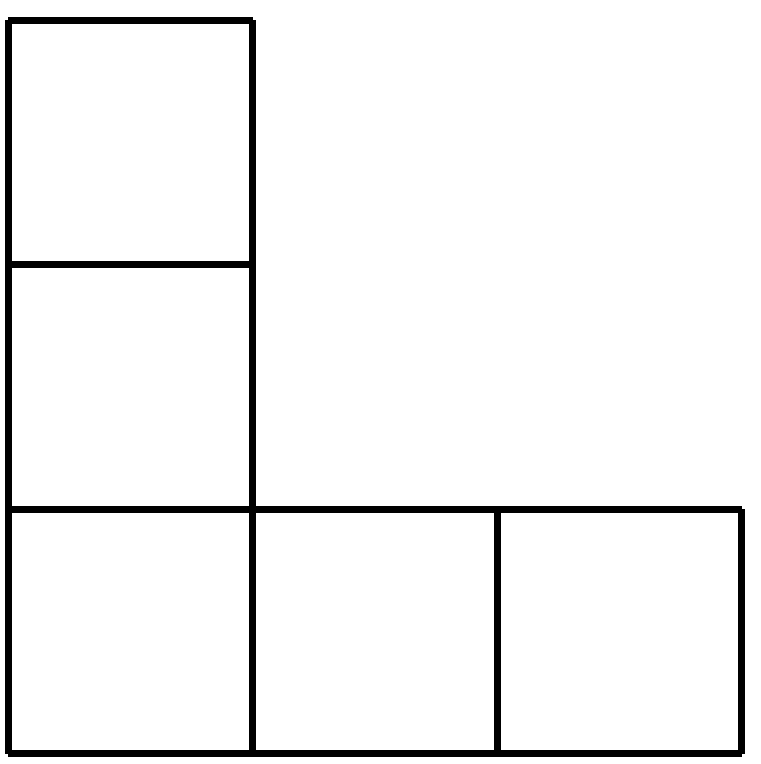}} &
\quad
\parbox{.17\textwidth}{\includegraphics[width=.20\textwidth]{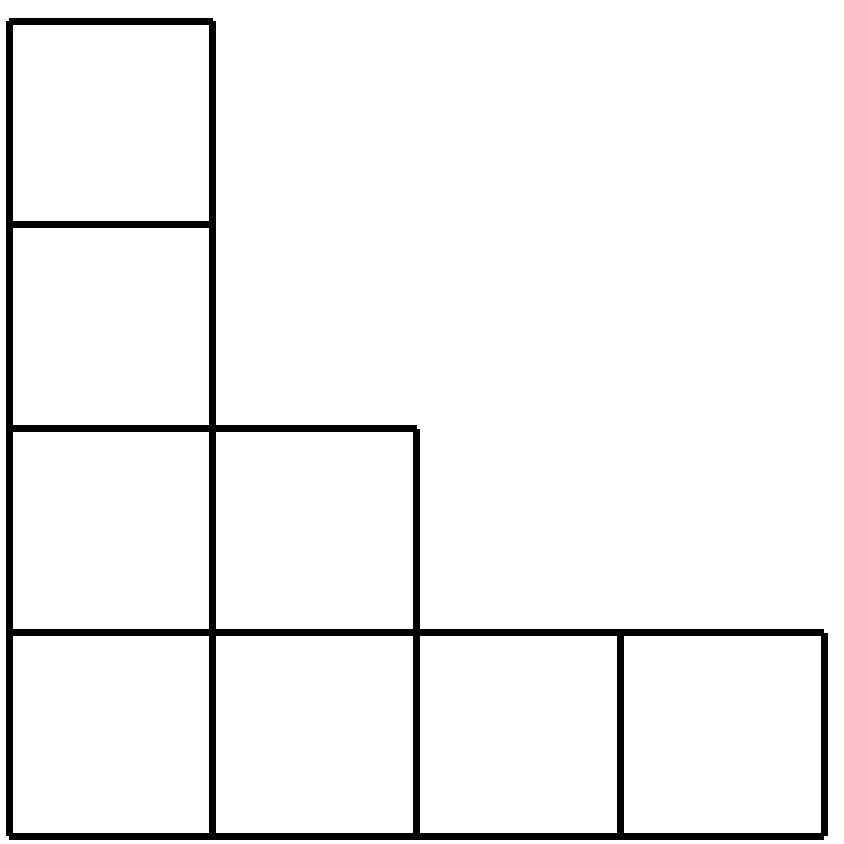}} 
\\
\\
\parbox{.17\textwidth}{\includegraphics[width=.20\textwidth]{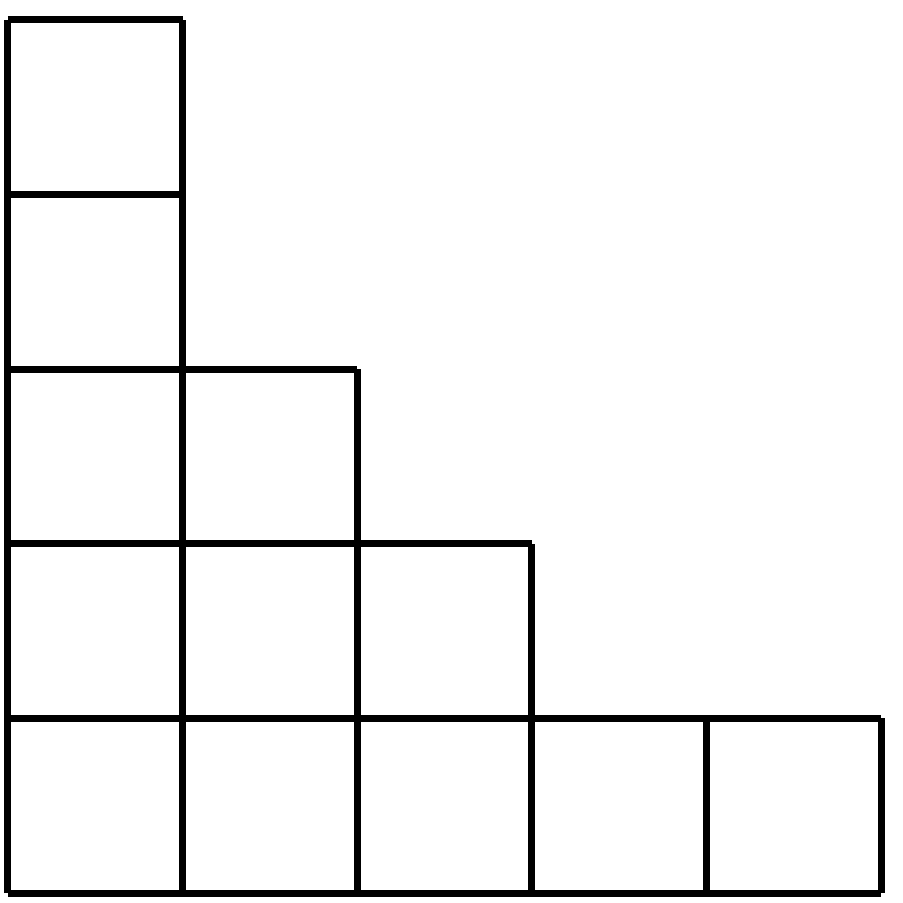}} &
\quad
\parbox{.17\textwidth}{\includegraphics[width=.20\textwidth]{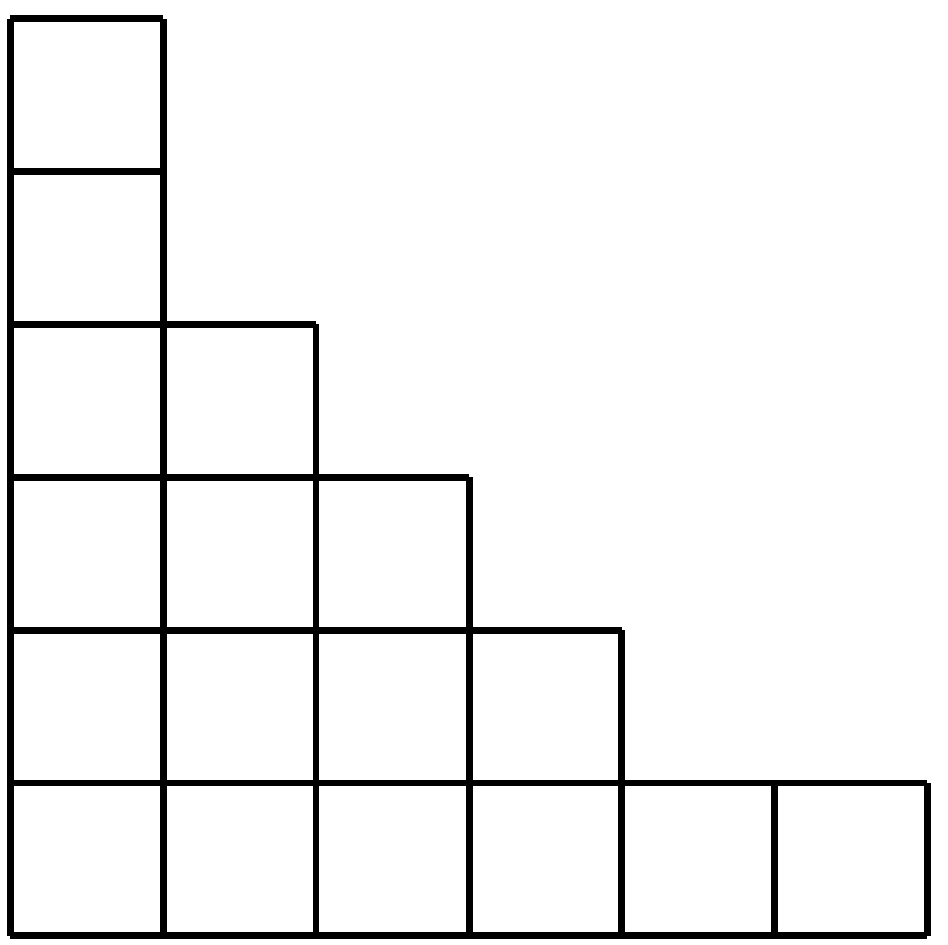}} &
\quad
\parbox{.17\textwidth}{\includegraphics[width=.20\textwidth]{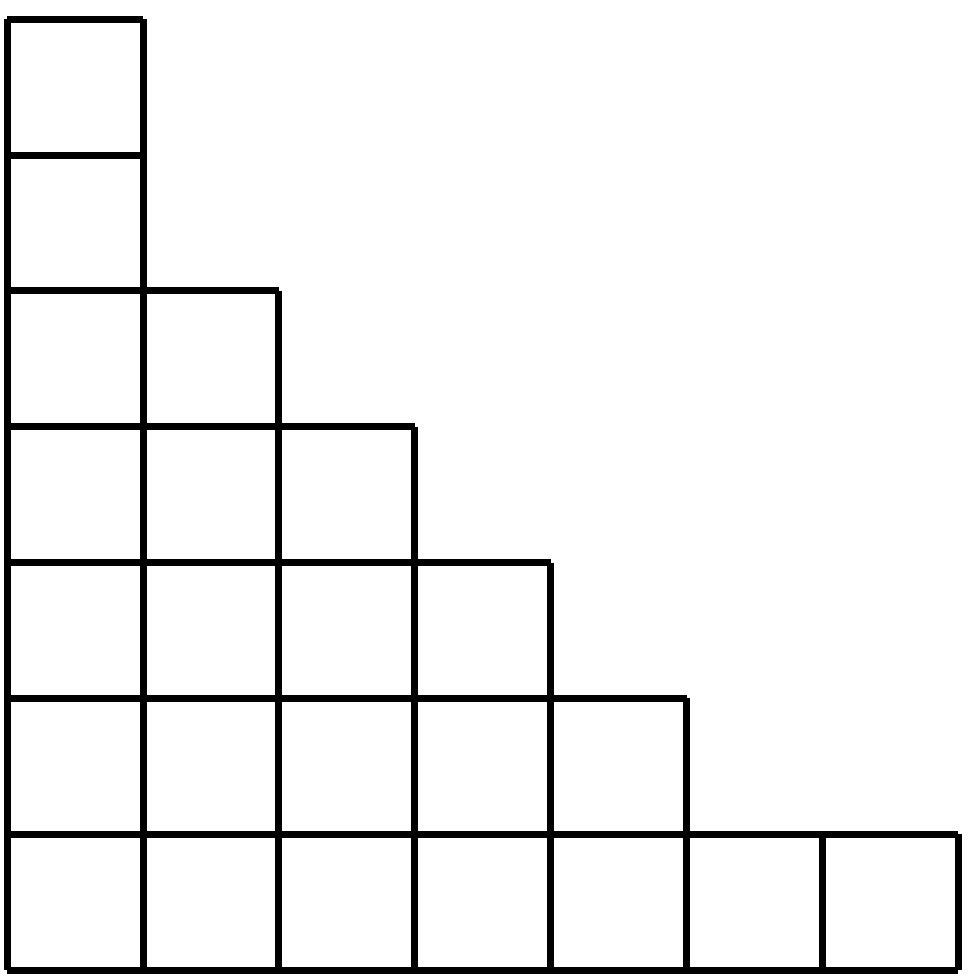}} 
\end{tabular}
}
\caption{For $n=2$, the geometry of the lower set $S_r$ is shown for $r=2,3,\ldots,7$.  Treating each figure as a set of unit squares with the lower left corner at the origin  in ${\RR}^2$, the corners of each square indicate the points of ${\NN}^2_0$ that belong to $S_r$.}
\label{fig:S2D}
\end{figure}

\begin{figure}
\centering{
\begin{tabular}{ccc}
\parbox{.26\textwidth}{\includegraphics[width=.26\textwidth]{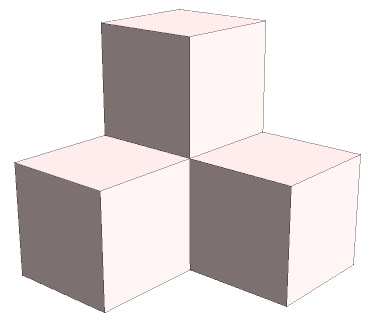}} &
\parbox{.26\textwidth}{\includegraphics[width=.26\textwidth]{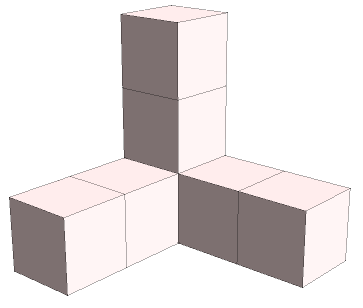}} &
\parbox{.26\textwidth}{\includegraphics[width=.26\textwidth]{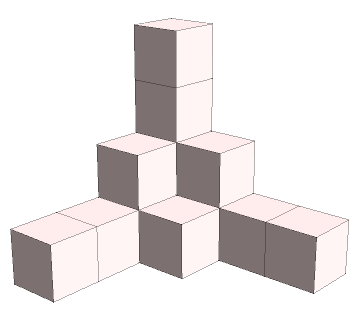}} 
\\
\parbox{.26\textwidth}{\includegraphics[width=.26\textwidth]{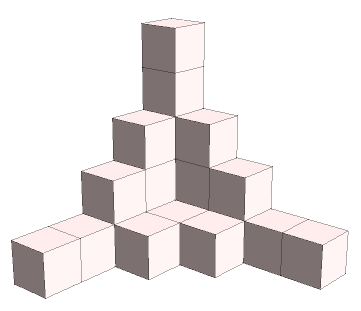}} &
\parbox{.26\textwidth}{\includegraphics[width=.26\textwidth]{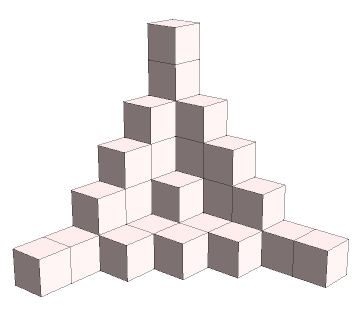}} &
\parbox{.26\textwidth}{\includegraphics[width=.26\textwidth]{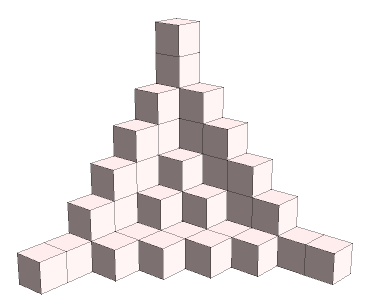}} 
\end{tabular}
}
\caption{For $n=3$, the geometry of the lower set $S_r$ is shown for $r=2,3,\ldots,7$.  Treating each figure as a set of unit cubes based at the origin and viewed from first octant in ${\RR}^3$, the corners of each cube indicate the points of ${\NN}^3_0$ that belong to $S_r$.}
\label{fig:S3D}
\end{figure}

We use this partition to compute the dimension of $S_r$ and confirm
that it agrees with the dimension of $\cS_r$ given in~\cite{AA}.
Fix $\beta\in\{0,1,2\}^n$ and let $d=\dim f_\beta$.
Letting $\NN_2$ denote natural numbers $\geq 2$, we see that
\[ \# S_{r,\beta} =  \#\{ \alpha \in \NN_2^d : |\alpha| \le r \} 
  = \#\{ \alpha \in \NN_0^d : |\alpha| \le r-2d \}. \]
Therefore,
\begin{equation}
\label{eq:crs-def}
 \# S_{r,\beta} = \begin{cases} \binom{r-d}{d}, & r \ge 2d; \\
                           0, & \hbox{otherwise}.
     \end{cases}
\end{equation}
Using (\ref{eq:numfaces}), we thus find
\[ \# S_r =
     \sum_{d=0}^n 2^{n-d} \binom{n}{d} \# S_{r,\beta} =
    \sum_{d=0}^{\min(n,\lfloor r/2 \rfloor)}
        2^{n-d} \binom{n}{d} \binom{r-d}{d} , \]
which is the formula for $\dim\cS_r$ in~\cite[Equation (2.1)]{AA}.
A table of values of $\dim \cS_r$ for small values of
$n$ and $r$ is given in \cite{AA}.
Figures~\ref{fig:S2D} and~\ref{fig:S3D} show the set $S_r$
for $r=2,3,\ldots,7$ in 2-D and 3-D respectively.

\section{Basis functions}\label{sec:nodal}

We now apply Theorem~\ref{thm:NewtonH} to the lower set $L = S_r$
to construct a nodal basis for $\cS_r(I^n)$ for arbitrary $r,n\geq 1$.
To do this, we choose the grid coordinates
$x_{j,k}$, $j=1,\ldots,n$, $k = 0,\ldots,r$,
in a manner that respects the indexing of the faces of $I^n$.
Suppose that for $j=1,\ldots,n$,
$$ x_{j,0} = -1 \quad \hbox{and} \qquad  x_{j,1} = 1, $$
and
$$ x_{j,k} \in (-1,1), \quad k=2,\ldots,r, $$
(not-necessarily distinct). Then
for each $\beta \in \{0,1,2\}^n$,
$x_\alpha \in f_\beta$ if and only if $\alpha \in S_{r,\beta}$.

Suppose further that the grid coordinates
$x_{j,k}$, $k=2,\ldots,r$, are distinct.
In this case the interpolation conditions of
Theorem~\ref{thm:NewtonH} are of Lagrange type:
\begin{equation}\label{eq:agreeslag}
 p(x_\alpha) = u(x_\alpha), \qquad \alpha \in S_r,
\end{equation}
giving the basis $\{\phi_{\alpha} : \alpha \in S_r \}$
for $\cS_r(I^n)$ defined by
$$ \phi_{\alpha}(x_{\alpha'}) = \delta_{\alpha,\alpha'},
  \qquad \alpha, \alpha' \in S_r. $$ 

\begin{figure}
\centering
\sbox{\strutbox}{\rule{0pt}{15pt}}           
\begin{tabular}[\textwidth]{@{\extracolsep{\fill}} ccc}
\parbox{.25\textwidth}
{\includegraphics[width=.25\textwidth]{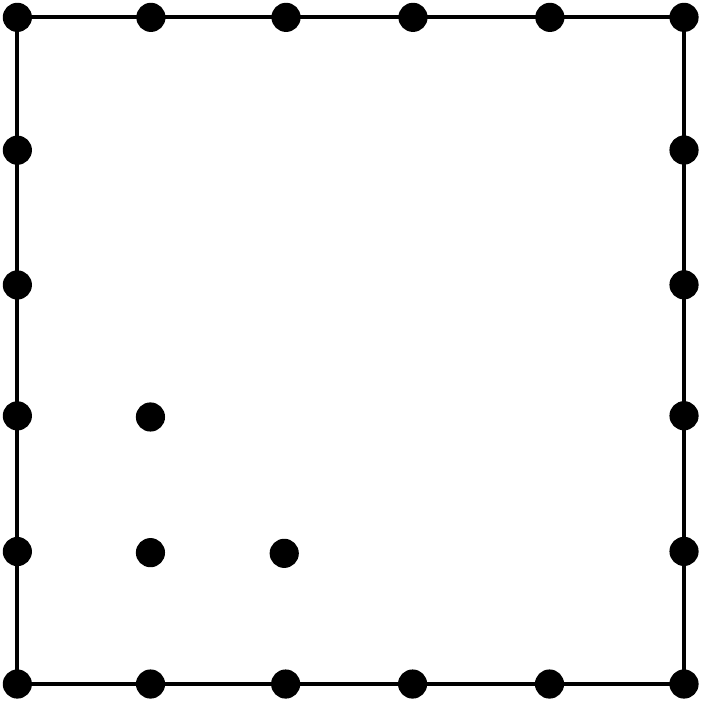}} & 
\parbox{.25\textwidth}
{\includegraphics[width=.25\textwidth]{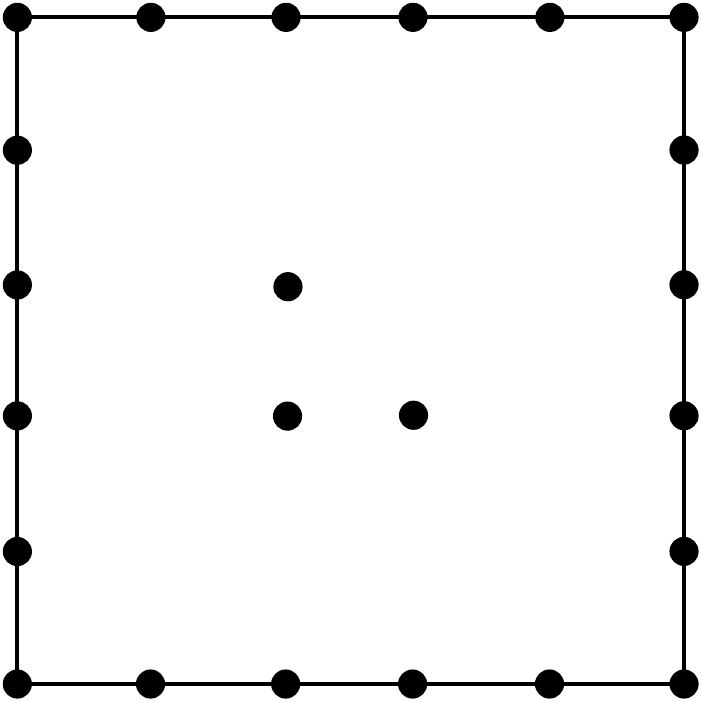}} & 
\parbox{.25\textwidth}
{\includegraphics[width=.25\textwidth]{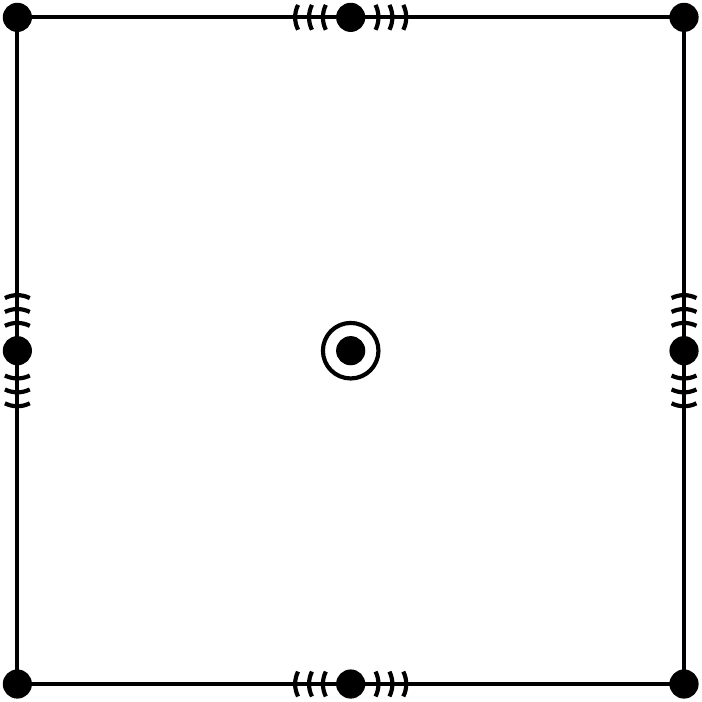}} \\
$(a)$ & $(b)$ & $(c)$
\end{tabular}
\caption{Three choices of $x_{j,k}$, $2 \le k \le 5$,
for $n=2$, $r=5$.
The first two choices are Lagrange-like while the third is Hermite-like.
}
\label{fig:partitions}
\end{figure}

\begin{figure}
\centering
\sbox{\strutbox}{\rule{0pt}{0pt}}           
\begin{tabular}[.64\textwidth]{@{\extracolsep{\fill}} cc}
\parbox{.32\textwidth}{\includegraphics[width=.32\textwidth]{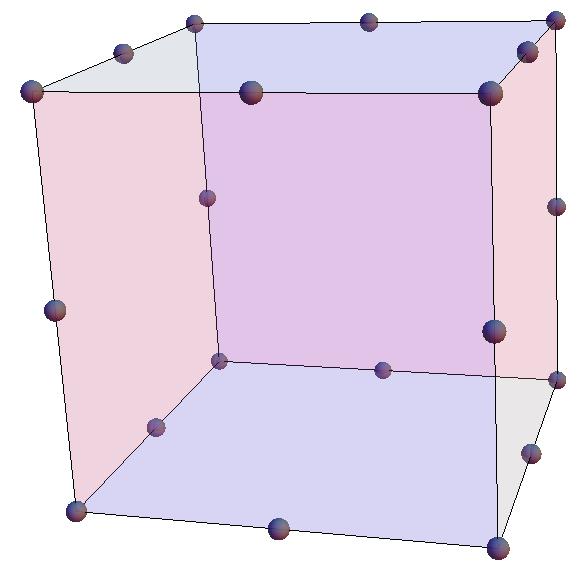}} & 
\parbox{.32\textwidth}{\includegraphics[width=.32\textwidth]{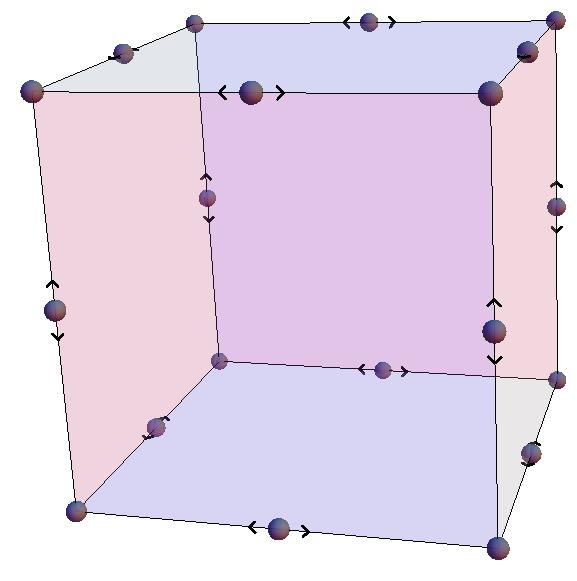}} \\
\parbox{.32\textwidth}{\includegraphics[width=.32\textwidth]{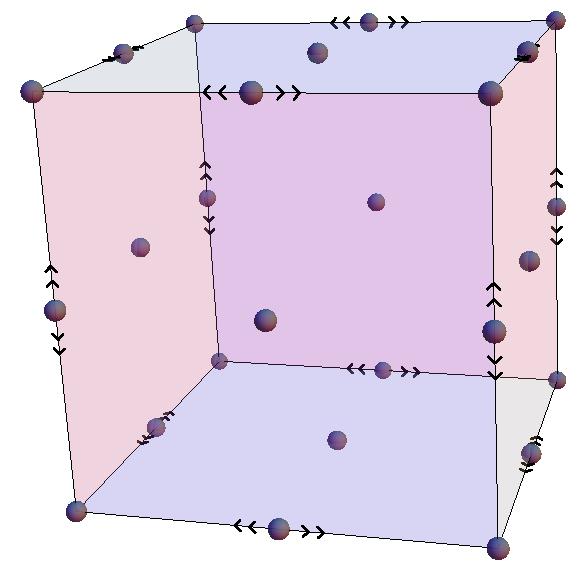}} &
\parbox{.32\textwidth}{\includegraphics[width=.32\textwidth]{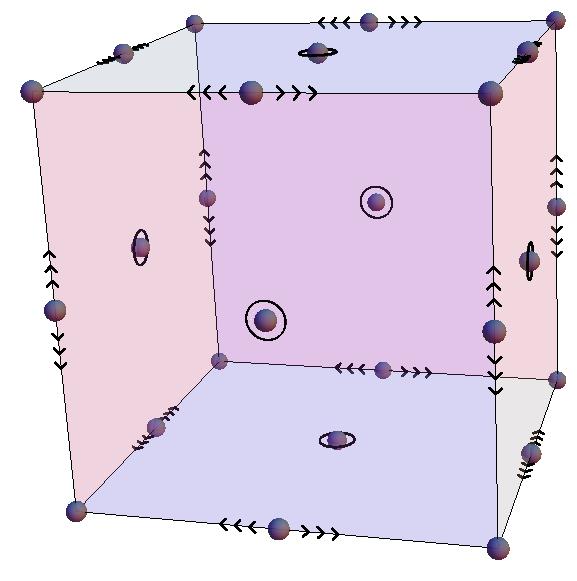}} \\ 
\end{tabular}
\caption{Hermite-like interpolation conditions in 3-D for $r=2,3,4,5$.
A dot indicates that a basis function will interpolate the value of the function at that location.
A dot on an edge enclosed by $\ell$ sets of parentheses indicates that basis functions will interpolate each partial derivative along the edge at the location of the dot, up to order $\ell$.
A dot in the interior enclosed by $\ell$ circles indicates that basis functions will interpolate all partial derivatives at the location of dot, up to total order $\ell$.}
\label{fig:IC3D}
\end{figure}

We consider two choices of such distinct coordinates.
The first choice is to distribute them uniformly in $I$
in increasing order:
\begin{equation}
\label{eq:xj-ex-a}
x_{j,k} = -1 + \frac{2(k-1)}{r},  \qquad k=2,\ldots,r,
\end{equation}
as illustrated in Figure~\ref{fig:partitions}a.
This configuration of nodes is, however, only symmetric for $r \le 3$.
Next, to obtain a more symmetric configuration,
we re-order the interior grid coordinates
in such a way that they are closer to the middle of $I$:
\begin{align}\label{eq:xj-ex-b}
x_{j,r-2s} &= 1 - \frac{2(s+1)}{r}, 
           \quad s=0,1,2,\ldots, \lfloor (r-2) / 2 \rfloor, \cr
x_{j,r-2s-1} &= -1 + \frac{2(s+1)}{r},
           \quad s=0,1,2,\ldots, \lfloor (r-3) / 2 \rfloor.
\end{align}
as illustrated in Figure~\ref{fig:partitions}b.
This yields a symmetric configuration for $r \le 4$,
but not for $r \ge 5$.
This lack of symmetry motivates the third choice of letting
all interior grid coordinates coalesce
to the midpoint of $I$, i.e.,
\begin{equation}
\label{eq:xj-ex-c}
x_{j,k} = 0, \qquad k=2,\ldots,r,
\end{equation}
as indicated in Figure~\ref{fig:partitions}c.
This gives interpolation conditions that are symmetric
for all $n$ and $r$, but the trade-off is that these conditions are
now of Hermite type rather than Lagrange.
In this Hermite case, all the points $x_\alpha$ in the
face $f_\beta$ are equal to the mipoint of that face, which we
denote by $y_\beta$.
The interpolation conditions of
Theorem~\ref{thm:NewtonH} then become
\begin{equation}
\label{eq:agreesherm}
 D^\rho p(y_\beta) = D^\rho u(y_\beta), \qquad
   \beta\in\{0,1,2\}^n, \quad \rho\in K_{r,\beta},
\end{equation}
where
\begin{equation}
K_{r,\beta} :=   \{\rho \in \NN_0^n : |\rho| \leq r - 2d~\text{with}~\rho_j
= 0~\text{if}~\beta_j<2\}.
\end{equation}
Thus $p$ can be expressed as
\[ p(x) = \sum_{\beta \in \{0,1,2\}^n} \sum_{\rho \in K_{r,\beta}}
  D^\rho u(y_\beta) \phi_{\beta,\rho}(x), \]
where
\[ \{\phi_{\beta,\rho} : \beta \in \{0,1,2\}^n, \rho \in K_{r,\beta} \} \]
is a basis for $\cS_r$ defined by
\[ D^{\rho'} \phi_{\beta,\rho}(y_{\beta'}) =
 \delta_{\beta,\beta'} \delta_{\rho,\rho'},\;\;\text{for any}\;\; 
   \beta' \in \{0,1,2\}^n, \;\rho' \in K_{r,\beta}. \]
Figure~\ref{fig:IC3D} illustrates
these interpolation conditions for $r=2$ through $r=5$
in the case $n=3$.



\section{Tensor-product formula}\label{sec:tp}

In this last section we explain how the interpolant
can be expressed as a linear combination of
tensor-product interpolants over various rectangular
subgrids of the overall grid.
This applies also to the basis functions and so gives
a simple method of evaluating these functions and their
derivatives.
To do this we apply the formula recently obtained in \cite{DF}.
Suppose again that $L \subset \NN_0^n$ is any lower set
as in Section~\ref{sec:lower} and
consider the interpolant $p$ to $u$ in
Theorem~\ref{thm:NewtonH}. 
For any $\alpha \in L$ define the rectangular block
$$ B_\alpha = \{\mu \in \NN_0^n : \mu \le \alpha \} $$
and let $p_\alpha \in P(B_\alpha)$ denote the
tensor-product interpolant to $u$ satisfying the
interpolation conditions (\ref{eq:agrees})
for $\mu \in B_\alpha$.
Further, let
$\chi(L) : \NN_0^n \to \{0,1\}$ be the characteristic
function defined by
$$ \chi(L)(\alpha) = \begin{cases} 1 & \hbox{if $\alpha \in L$;} \\
                                  0 & \hbox{otherwise.} \\
  \end{cases} $$
It was shown in \cite{DF} that
\begin{equation}\label{eq:pI}
 p = \sum_{\alpha \in L} c_\alpha p_\alpha,
\end{equation}
where
\begin{equation}\label{eq:csimple}
 c_\alpha = \sum_{\epsilon \in \{0,1\}^n} (-1)^{|\epsilon|} 
                \chi(L)(\alpha + \epsilon), \qquad \alpha \in L.
\end{equation}
\noindent This in turn gives a formula for each basis function
$\phi_\beta \in P(L)$, i.e.,
\begin{equation}
\label{eq:nodalformula}
 \phi_\beta(x) = \sum_{\substack{\alpha \in L \\ \alpha \ge \beta}}
                   c_\alpha \phi_{\beta,\alpha},
\end{equation}
where $\phi_{\beta,\alpha} \in P(B_\alpha)$ denotes the tensor-product
basis function associated with the index $\beta$,
defined by
$$ \lambda_{\beta'} \phi_{\beta,\alpha} = \delta_{\beta,\beta'},
     \qquad \beta \in B_\alpha. $$
For a general lower set $L$, many of the integer coefficients
$c_\alpha$ are zero, and so in order to apply
(\ref{eq:pI}) to evaluate $p$ we need to determine which of the
$c_\alpha$ are non-zero, and to find their values.
With $L=S_r$ we could do this in practice by
implementing the formula (\ref{eq:csimple}).
However, we will derive a specific formula for
the $c_\alpha$.
We call $\alpha \in L$ a \emph{boundary point} of $L$
if $\alpha + 1_n \not\in L$, where $1_n = (1,1,\ldots,1) \in \NN_0^n$.
Let $\partial L$ denote the set of boundary points of $L$.
As observed in \cite{DF}, if $\alpha$ is not a boundary point then
$c_\alpha = 0$.

Consider now the formula (\ref{eq:pI}) when $L=S_r$.
Note that $|\alpha|'$ is a symmetric function of $\alpha$: it is unchanged
if we swap $\alpha_j$ and $\alpha_i$ for $i \ne j$.
It follows that $\chi(S_r)(\alpha)$ is also symmetric in
$\alpha$, and therefore $c_\alpha$ is also symmetric in $\alpha$.
We can thus determine the boundary points $\alpha \in \partial S_r$ and their
coefficients $c_\alpha$ according to how many zeros and
ones $\alpha$ contains.
For any $\alpha \in \NN_0^n$ let
$m_i(\alpha)$ denote the multiplicity of the integer
$i \ge 0$ in $(\alpha_1,\ldots,\alpha_n)$, i.e.,
$$ m_i(\alpha) = \# \{\alpha_j = i\}. $$

\begin{lemma}\label{lem:alphazero}
If $\alpha \in \partial S_r$ and $m_0(\alpha) \ge 1$ then $c_\alpha = 0$.
\end{lemma}

\begin{proof}
By the symmetry of $c_\alpha$ we may assume that
$\alpha_1 = 0$, and from (\ref{eq:csimple}) we can express $c_\alpha$ as
$$
 c_\alpha = \sum_{\epsilon \in \{0\} \times \{0,1\}^{n-1}} (-1)^{|\epsilon|} 
                \big( \chi(S_r)(\alpha + \epsilon)
     - \chi(S_r)(\alpha + e_1 + \epsilon) \big),
$$
where $e_1 = (1,0,\ldots,0) \in \NN_0^n$.
Since $\alpha_1 = 0$, both
$$ (\alpha + \epsilon)_1 \le 1 \qquad \hbox{and} \qquad
   (\alpha + e_1 + \epsilon)_1 \le 1, $$
and so
$$ |\alpha+\epsilon|' = |\alpha+e_1+\epsilon|', $$ 
and therefore
$$ \chi(S_r)(\alpha + \epsilon) = \chi(S_r)(\alpha + e_1 + \epsilon), $$
and so $c_\alpha =0$.
\end{proof}

In view of Lemma~\ref{lem:alphazero},
we need only consider points
$\alpha \in \partial S_r \cap \NN_1^n$.

\begin{lemma}\label{lem:boundarypoints}
Let $\alpha \in S_r \cap \NN_1^n$ and $m_1 = m_1(\alpha)$.
Then $\alpha \in \partial S_r$
if and only if
$$ |\alpha|' > r - (n + m_1). $$
\end{lemma}

\begin{proof}
By the definition of $S_r$,
$\alpha \in \partial S_r$ if and only if $|\alpha+1_n|' > r$.
Since $\alpha \in \NN_1^n$,
$$ \# \{\alpha_j \ge 2\} =  n - m_1, $$
and we find
$$ |\alpha+1_n|' = 2m_1 + |\alpha|' + (n-m_1)
   = |\alpha|' + n + m_1, $$
which proves the result.
\end{proof}

In view of Lemma~\ref{lem:boundarypoints},
we need only consider points $\alpha \in \NN_1^n$ such that
\begin{equation}\label{eq:alphak}
 |\alpha|' = r-k, \qquad k=0,1,\ldots, n + m - 1,
\end{equation}
where $m = m_1(\alpha)$.

\begin{theorem}\label{thm:calpha}
Let $\alpha \in \NN_1^n$ be as in (\ref{eq:alphak}).
If $m < n$ then
\begin{equation}\label{eq:cmk}
   c_\alpha = c_{m,k} :=
  \sum_{i=0}^{m} 
  (-1)^{k+i}
  \binom{m}{i}
  \binom{n-m-1}{k-2i},
\end{equation}
with the convention that
$\binom{l}{j} = 0$ if $j < 0$ or $j > l$.
\end{theorem}

\begin{proof}
Let $\epsilon \in \{0,1\}^n$, and let
$$ i_1 = \# \{j : \hbox{$\epsilon_j = 1$ and $\alpha_j = 1$}\}, $$
$$ i_2 = \# \{j : \hbox{$\epsilon_j = 1$ and $\alpha_j \ge 2$}\}. $$
Then
$$ |\alpha+\epsilon|' = |\alpha|' + 2i_1 + i_2, $$
and so $\alpha+\epsilon \in S_r$ if and only if
$$ |\alpha|' + 2i_1 + i_2 \le r, $$
or, equivalently,
$$ 2i_1 + i_2 \le k. $$
Since the number of ways of choosing $i_1$ elements among $m$ is
$\binom{m}{i_1}$, and
the number of ways of choosing $i_2$ elements among $n-m$ is
$\binom{n-m}{i_2}$ the sum in (\ref{eq:csimple})
reduces to
$$
   c_\alpha =
  \sum_{i_1=0}^{m} 
  \sum_{i_2=0}^{k-2i_1}
  (-1)^{i_1+i_2}
  \binom{m}{i_1}
  \binom{n-m}{i_2}.
$$
Since
$$
  \sum_{i_2=0}^{k-2i_1}
  (-1)^{i_2}
  \binom{n-m}{i_2}
  =
  (-1)^{k-2i_1}
  \binom{n-m-1}{k-2i_1},
$$
we obtain (\ref{eq:cmk}).
\end{proof}

Table~\ref{tab:coeffs} shows
the values of the coefficients
$c_{m,k}$ for $n=1,2,3,4$.
\begin{table}[h]
\begin{center}
{\small
\begin{tabular}{r|r|rrrrrrr}
  & & & & & $k$ \\
$n$ & $m$ & $0$ & $1$ & $2$ & $3$ & $4$ & $5$ & $6$ \\
\hline
& & & \\
\hline
$1$
& 0 &  1 \\
\hline
& & & \\
\hline
$2$
& 0 &  1 & -1 \\
& 1 &  1 &  0 & -1 \\
\hline
& & & \\
\hline
$3$
& 0 &  1 & -2 &  1 \\
& 1 &  1 & -1 & -1 &  1 \\
& 2 &  1 &  0 & -2 &  0 &  1 \\
\hline
& & & \\
\hline
$4$
& 0 &  1 & -3 &  3 & -1 \\
& 1 &  1 & -2 &  0 &  2 & -1 \\
& 2 &  1 & -1 & -2 &  2 &  1 & -1 \\
& 3 &  1 &  0 & -3 &  0 &  3 &  0 & -1
\end{tabular}
}
\end{center}
\caption{Coefficients $c_{m,k}$ for $n=1,2,3,4$.}
\label{tab:coeffs}
\end{table}
Finally, we need to consider the possibility that $m=n$ in
(\ref{eq:alphak}), in which case the formula (\ref{eq:cmk}) is
no longer valid, and we must treat this situation separately.
In this case $\alpha = 1_n$ and we can again
find $c_\alpha$ from (\ref{eq:csimple}).
Since 
$$ |1_n + 1_n|' = |2_n|' = 2n, $$
we see that $1_n \in \partial S_r$ if and only if $r < 2n$.

\begin{theorem}\label{thm:c1}
Suppose that $r < 2n$. Then
\begin{equation}\label{eq:c1}
   c_{1_n} = (-1)^{\lfloor r/2 \rfloor} 
  \binom{n-1}{\lfloor r/2 \rfloor}.
\end{equation}
\end{theorem}

\begin{proof}
For any $\epsilon \in \{0,1\}^n$,
$$ |1_n + \epsilon|' = 2|\epsilon|, $$ and so
(\ref{eq:csimple})
gives
$$
 c_{1_n}
  = \sum_{\substack{\epsilon \in \{0,1\}^n \\ 2|\epsilon| \le r}}
   (-1)^{|\epsilon|} 
  =
  \sum_{i=0}^{\lfloor r/2 \rfloor}
  \binom{n}{i} (-1)^i
$$
which gives (\ref{eq:c1}).
\end{proof}

We now consider examples of the use
of Theorems~\ref{thm:calpha} and~\ref{thm:c1},
and let $p_r$ denote the interpolant $p$ in
Theorem~\ref{thm:NewtonH} when $L=S_r$.

\subsection{2-D case}
\begin{figure}[t]
\centering{
\includegraphics[width=0.20\linewidth]{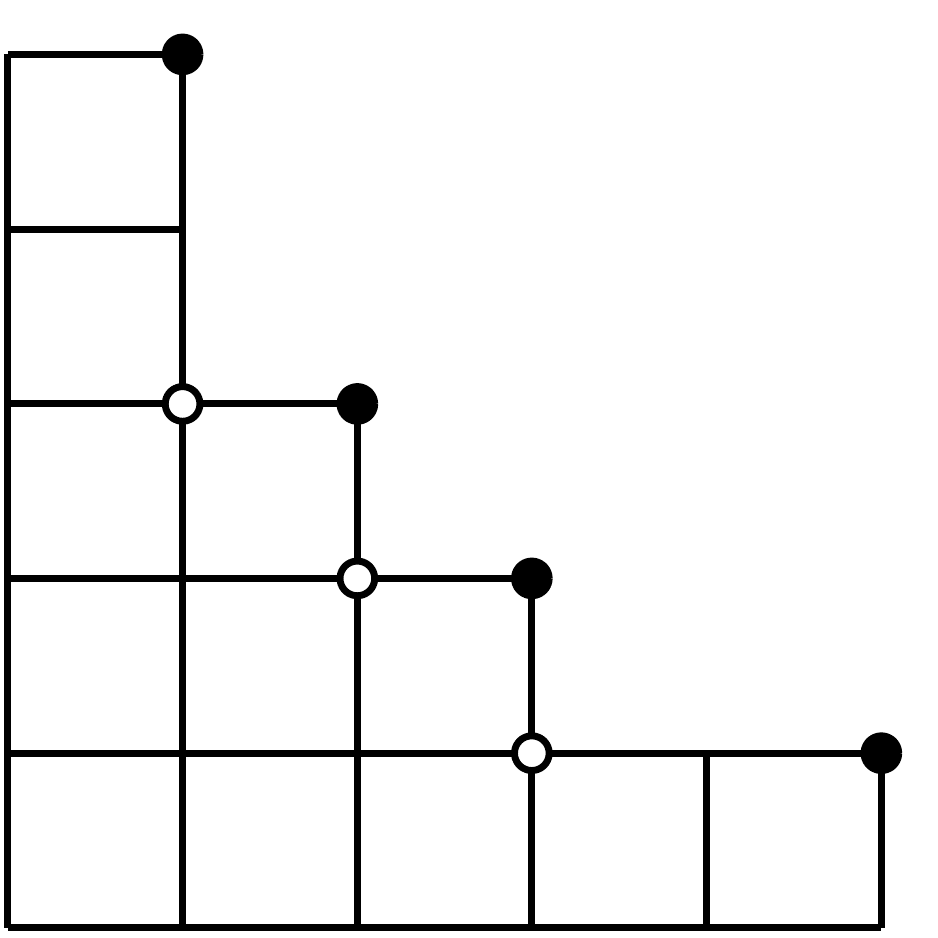}
}
\par
\caption{The geometry of $S_5$ for $n=2$ is shown (see Figure~\ref{fig:S2D}) with an indication of which blocks within the set contribute to the representation of the serendipity interpolant $p_5$ as a linear combination of tensor product interpolants.
A block with a filled dot in the upper right corner contributes with coefficient +1 while a block with an empty dot in the upper right corner contributes with coefficient $-1$.}
\label{fig:S5TP}
\end{figure}

\begin{figure}[t]
\centering{
\begin{tabular}{cccccccccc}
  &
\parbox{.17\textwidth}
{\includegraphics[width=.17\textwidth]{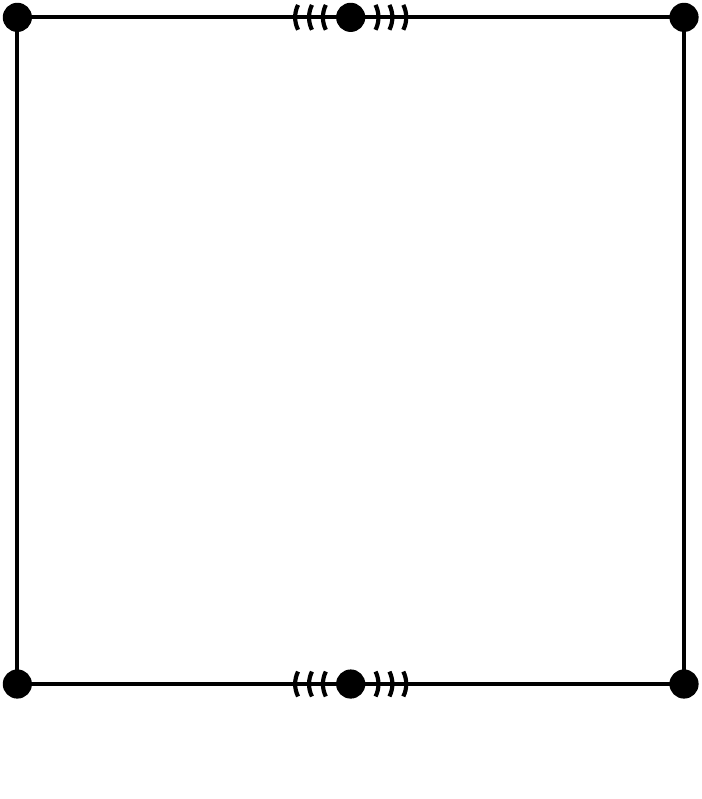}} &
$+$ &
\parbox{.17\textwidth}
{\includegraphics[width=.17\textwidth]{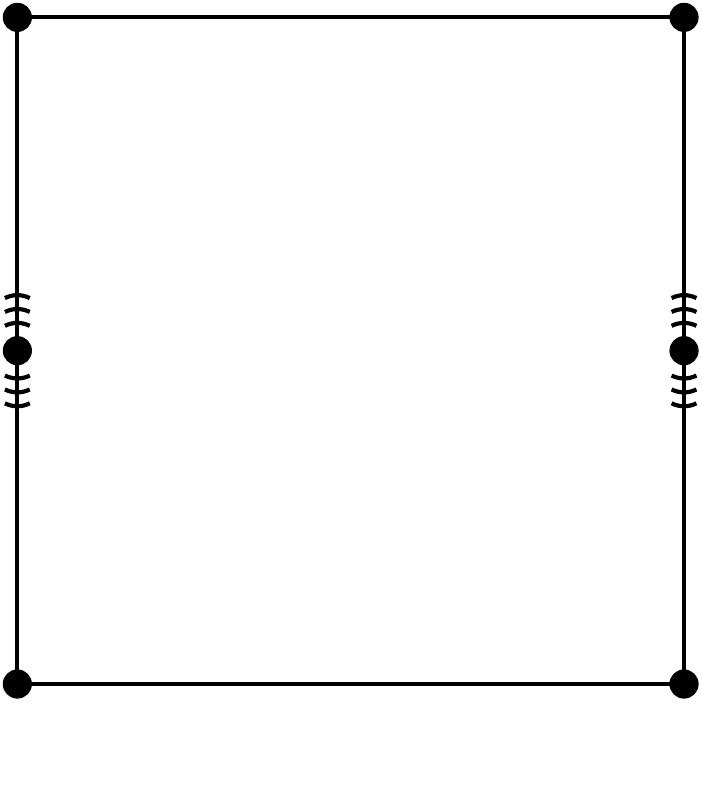}} &
$+$ &
\parbox{.17\textwidth}
{\includegraphics[width=.17\textwidth]{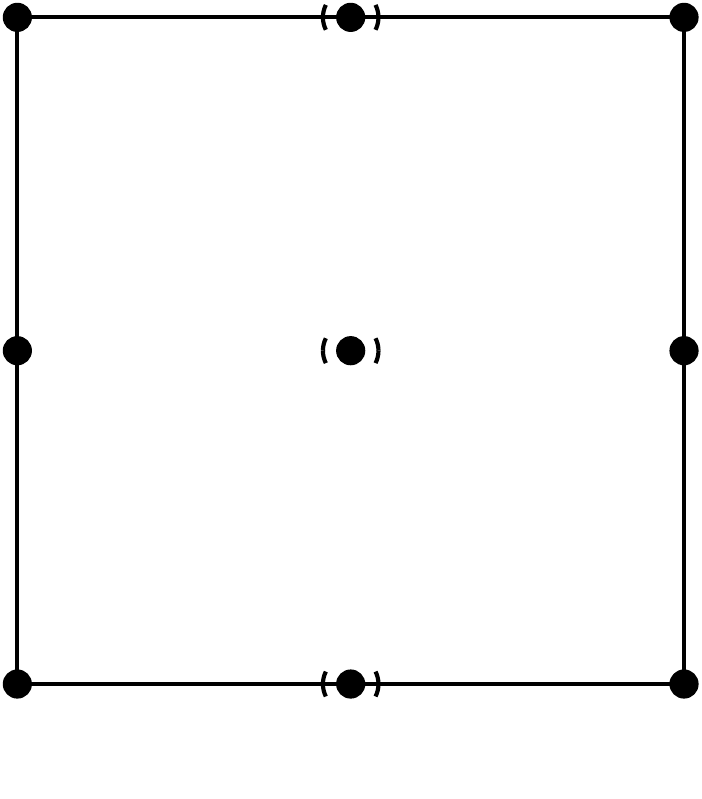}} & 
$+$ &
\parbox{.17\textwidth}
{\includegraphics[width=.17\textwidth]{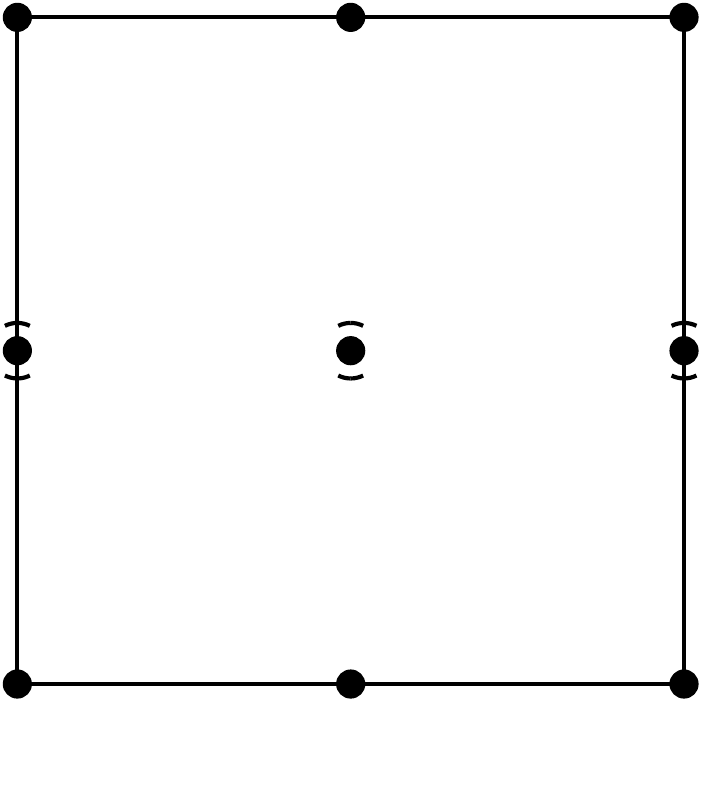}} & 
\\
\\
$-$ &
\parbox{.17\textwidth}
{\includegraphics[width=.17\textwidth]{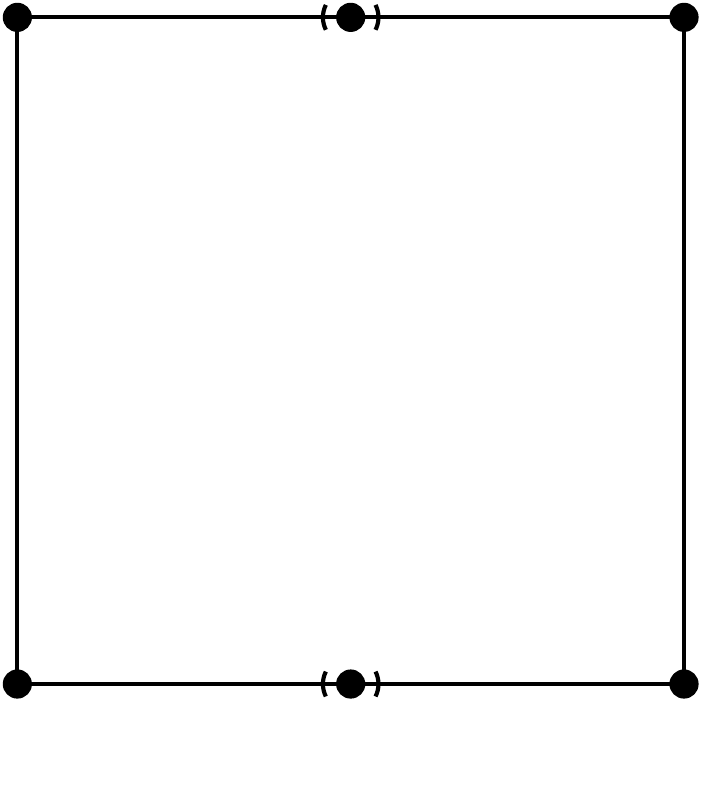}} &
$-$ &
\parbox{.17\textwidth}
{\includegraphics[width=.17\textwidth]{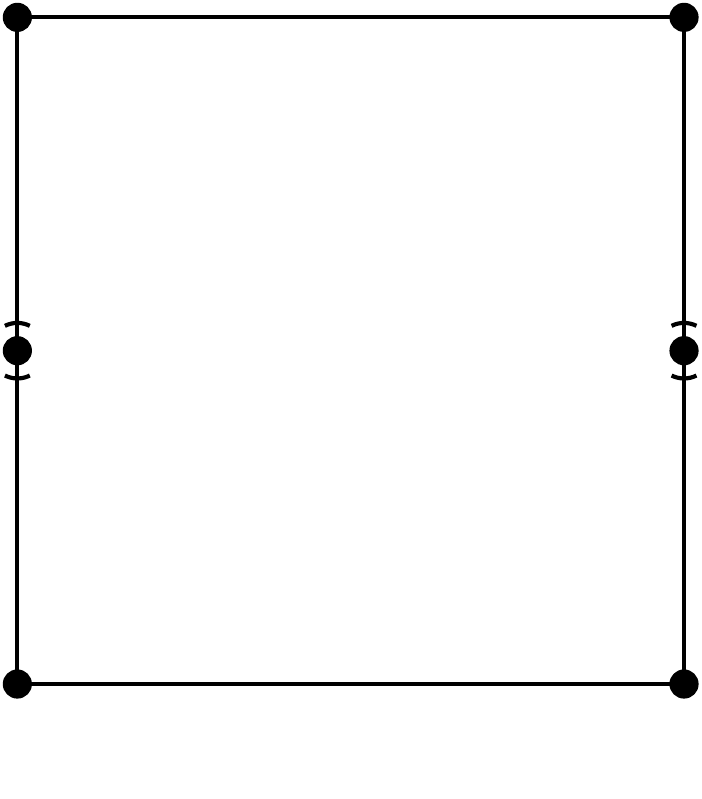}} &
$-$ &
\parbox{.17\textwidth}
{\includegraphics[width=.17\textwidth]{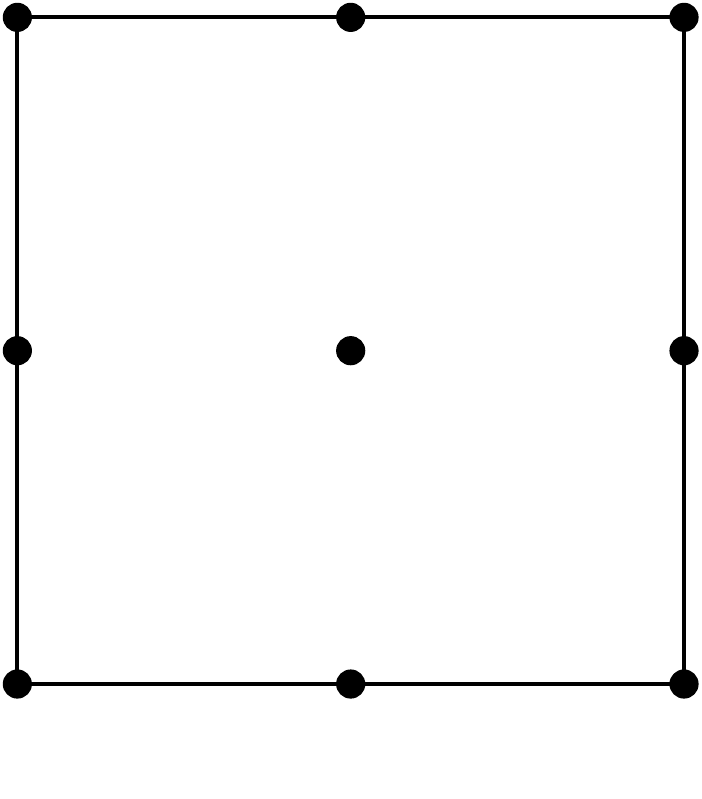}} &
$=$ &
\parbox{.17\textwidth}
{\includegraphics[width=.17\textwidth]{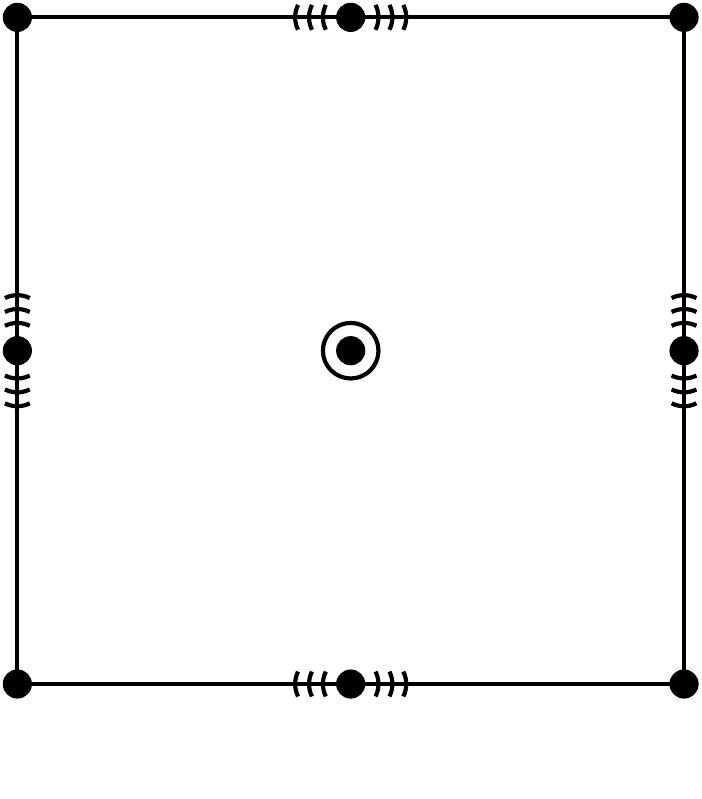}} 
\end{tabular}
}
\caption{A visual depiction of the the formula for $p_5$ in the Hermite case.}
\label{fig:p5herm}
\end{figure}

For $n=2$, Theorems~\ref{thm:calpha} and~\ref{thm:c1} give
\begin{align*}
 p_1 &= p_{11}, \cr
 p_2 &= p_{21} + p_{12} - p_{11}, \cr
 p_3 &= p_{31} + p_{13} - p_{11}, \cr
 p_4 &= p_{41} + p_{14} + p_{22} - (p_{21} + p_{12}), \cr
 p_5 & = p_{51} + p_{15} + p_{32} + p_{23} - (p_{31} + p_{31} + p_{22}).
\end{align*}
Figure~\ref{fig:S5TP} shows the polynomials in $S_5$
in the formula for $p_5$,
with black if $c_\alpha = 1$ and white if $c_\alpha = -1$.
%
%
Figure~\ref{fig:p5herm} depicts the polynomials in
the same formula, based on the Hermite interpolation conditions
(\ref{eq:xj-ex-c}).

\subsection{3-D case}
For $n=3$, to simplify the formulas let
$$ q_\alpha := \sum_{\alpha' \in \pi(\alpha)} p_{\alpha'}, $$
with $\pi(\alpha)$ denoting all permutations of
$\alpha = (\alpha_1,\alpha_2,\alpha_3)$, so that, for example,
\begin{align*}
 q_{111} &: = p_{111}, \cr
 q_{112} &: = p_{112} + p_{121} + p_{211}, \cr
 q_{123} &: = p_{123} + p_{132} + p_{213}
            + p_{231} + p_{312} + p_{323},
\end{align*}
etc. Then Theorems~\ref{thm:calpha} and~\ref{thm:c1} give
\begin{align*}
 p_1 &= q_{111}, \cr
 p_2 &= q_{112} - 2 q_{111}, \cr
 p_3 &= q_{113} - 2 q_{111}, \cr
 p_4 &= q_{122} + (q_{114}- 2 q_{112}) + q_{111}, \cr
 p_5 &= (q_{123} - q_{122}) + (q_{115}- 2 q_{113}) + q_{111}.
\end{align*}
We note that Delvos \cite{D} found a nodal basis for $p_4$, $n=3$,
using his method of `Boolean interpolation.'
That method is not, however, general enough to give
the formulas for $p_r$ with $r \ge 5$, $n = 3$.
Now that we have provided a generalized approach to defining nodal bases for serendipity elements, it remains to be studied whether certain arrangements of the grid coordinates $x_{j,k}$ provide advantages in specific application contexts.
Suitable pre-conditioners associated to these bases may also be needed.



\bibliographystyle{spmpsci}
\bibliography{srdp}

\end{document}